\documentclass[12pt]{article}

\usepackage{makeidx}
\usepackage[centertags]{amsmath}
\usepackage{amsthm}
\usepackage{newlfont}
\usepackage{amssymb}
\usepackage{amsmath}
\usepackage{amsbsy}
\usepackage{xcolor}
\usepackage{graphicx}
\usepackage{subfigure}
\usepackage{caption}
\usepackage{float}
\usepackage{rotating}
\usepackage{tablefootnote}
\usepackage{cases}
\usepackage{latexsym}
\usepackage{CJK}
\usepackage{epstopdf}
\usepackage{longtable}
\usepackage{supertabular}
\usepackage{etex}
\usepackage{amsfonts}
\usepackage{epsfig}
\usepackage{array}
\usepackage{verbatim}
\usepackage{bm}
\usepackage{manyfoot}
\usepackage{footmisc}
\usepackage{ctable}
\usepackage{booktabs}
\usepackage{morefloats}
\usepackage{graphicx}
\usepackage{xcolor}
\usepackage{graphicx}
\usepackage{subfigure}
\usepackage{caption}
\usepackage{extarrows}
\usepackage{pdflscape}
\usepackage{enumerate}
\usepackage{natbib}
\usepackage[colorlinks, citecolor=blue,urlcolor=blue,pagebackref]{hyperref}

\def\proclaim#1{\par \bigskip\noindent {\bf #1}\bgroup\it\ }
\def\endproclaim{\egroup\par\bigskip}

\newcommand{\no}{\nonumber}

\newcommand{\ep}{\epsilon}
\newcommand{\la}{\label}
\newcommand{\be}{\begin{eqnarray}}
\newcommand{\ee}{\end{eqnarray}}
\newcommand{\bestar}{\begin{eqnarray*}}
\newcommand{\eestar}{\end{eqnarray*}}
\newcommand{\lam}{\lambda}

\newcommand{\al}{\alpha}
\newcommand{\de}{\delta}

\newcommand{\ga}{\gamma}

\newcommand{\si}{\sigma}

\newcommand{\ignore}[1]{}

\newcommand{\E}{\mathbb E}

\renewcommand{\(}{\Big(}
\renewcommand{\)}{\Big)}

\newbox\TempBox \newbox\TempBoxA

\def\text#1{\mbox{\rm #1}}

\def\underwiggle 1{
\ifmmode\setbox\TempBox=\hbox{$ 1$}\else\setbox\TempBox=\hbox{
1}\fi \setbox\TempBoxA=\hbox to \wd\TempBox{\hss\char'176\hss}
\rlap{\copy\TempBox}\smash{\lower9pt\hbox{\copy\TempBoxA}} }

\newtheorem{thm}{Theorem}[section]

\newtheorem{lem}{Lemma}[section]

\newtheorem{cor}{Corollary}[section]

\newtheorem{assump}{Assumption}[section]
\newtheorem{rem}{Remark}[section]

\newtheorem{exap}{Example}[section]

\begin{document}

\title{Limit theorems for a  class of  martingale arrays with applications in nonlinear cointegrating regression}
\author{ Qiying Wang \\
The University of Sydney \footnote{Address correspondence to: School of
Mathematics and Statistics, The University of Sydney, NSW 2006,
Australia. Email: qiying.wang@sydney.edu.au.} }

\maketitle

\baselineskip=0.7 true cm

\begin{abstract}

 This paper develops a new asymptotic theory for a broad class of martingales, establishing convergence to limiting distributions that involve a functional of stochastic integrals. The proposed limit theorem substantially extends existing martingale asymptotic theory by accommodating a wider class of  dependence structures. As a primary application, the theory is applied to nonlinear regression models with nonstationary time series, yielding a rigorous framework for asymptotic inference on nonlinear least square estimators.

\end{abstract}

\section{Introduction}

Let $\{\ep_{nj}, {\cal F}_{nj}\}_{j\ge 1, n\ge 1}$ be a    martingale difference array defined on a common probability space $(\Omega, {\cal F}, P)$  with $ \E\big(\ep_{nj}^2 | {\cal F}_{n, j-1} \big)=1, a.s.$.  For each $n\ge 1$, let $\{u_{nj}\}_{j\ge 1}$ and $\{x_{nj}\}_{j\ge 1}$ be arbitrary random  sequences  so that both are   adapted to ${\cal F}_n=\{{\cal F}_{n,j-1}\}_{j\ge 1}$ and, on $D_{{\mathbb R}^2}[0,\infty),$
\bestar
\big( \frac 1{\sqrt n }\sum_{j=1}^{[nt]} \ep_{nj},\ x_{n, [nt]}\big) &\Rightarrow& (B_t, X_t),
\eestar
where $B=\{B_t\}_{t\ge 0}$ is a  Brownian motion and $X=\{X_t\}_{t\ge 0}$ is  a continuous   process.  Define a  martingale array  by  $\{S_{nk}, {\cal F}_{nk}\}_{k\ge 2, n\ge 1}$, where
\bestar
S_{nk}=\frac 1{\sqrt n}\sum_{j=1}^k H(x_{nj}, u_{nj})\, \ep_{nj}, \quad k=2, 3, ...
\eestar
and $H$ is a continuous function of its components. In this paper, we concern with  the asymptotics of a statistic $S_n$ defined  by $S_n:=S_{nn}$, as $n\to\infty$ and their applications in nonlinear cointegrating regression.

The limit theory for martingales and/or statistics like $S_n$ derived from martingale arrays is a celebrated achievement in the development of Probability and Statistics. For various martingales with different structures, extensive articles have established a number of fundamental results on  central limit theorems, convergence to   a mixture of normal distributions and convergence to stochastic integrals.  The earlier contributions to these topics  were   summarized in excellent  books by \cite{HH1}, \cite{LS2} and \cite{JS1}.

 To the best of our knowledge, however, the asymptotic behavior of
 $S_n$ has not been investigated in the literature, because of the presence of both stationary
  (e.g.,$u_{nj}$) and non-stationary (e.g.,$x_{nj}$) variables.  First of all, when $u_{nj}$ is stationary,     no limitation exists for $H(x_{n, [nt]}, u_{n, [nt]})$ as $n\to \infty$  and therefore   the classical theory on the convergence of stochastic integrals (e.g., \cite{KP1}) is unlikely to  be used for such a martingale array $\{S_{nk}, {\cal F}_{nk}\}_{k\ge 2, n\ge 1}$.
On the other hand, by writing the conditional variance $V_n^2$ of $S_n$ by
$
V_n^2 = \frac {1}{ n}\, \sum_{j=1}^n H^2(x_{nj}, u_{nj}) ,
$
we may prove that $V_n^2\to_D \int_0^1 \widehat H^2(X_t) dt$ as  in Remark \ref {r2.2}, where  $\widehat H^2(x)$ is defined as in Assumption \ref {as2.3} below, but
 it is usually impossible to show  that $V_n^2\to_P  \int_0^1 \widehat H^2(X_t) dt$ under natural settings on $H(x, y)$ and $u_{nj}$\footnote{For a simple illustration, assume that $H(x, y)$ is a  continuous function on ${\mathbb R}^2$ and $x_{nj}=\frac 1{\sqrt n}\sum_{j=1}^{j}\ep_i$, where $\{\ep_j\}_{j\ge 1}$  is a sequence of iid random variables with $E\ep_1=0$ and $E\ep_1^2=1$. We further assume that  $u_{nj}=\ep_j$ and $\ep_{nj}=\ep_j$ so that
$
V_n^2=  \frac {1}{ n}\, \sum_{j=1}^n H^2(x_{nj}, \ep_{j}).
$
Note that $x_{n, [nt]}\Rightarrow B_t$ on $D[0,1]$, but the convergence in distribution can not be enhanced to the convergence in probability.
We may prove that $V_n^2\to_D \int_0^1 \hat H(B_t) dt$,  where $\hat H(x)= EH^2 (x, \ep_1)$, but it is incorrect to say that $V_n^2\to_P \int_0^1 \hat H(B_t) dt.$}. This claim suggests that  the classical martingale limit theorems are unable to establish the asymptotics of $S_n$. Indeed,   many  earlier results on  the convergence to  a mixture of normal distributions and on the convergence to  a  semimartingale    were established under the assumption  that $V_n^2 \to_P V^2$, where $V^2$ is an almost finite random variable  (e.g.,\cite{HH1}; \cite{JS1}).
More recently,  in investigating the convergence to a mixture of normal distributions for a martingale array,  \cite{J1},  \cite{W1,W2} and  \cite{RV1,RV2,RV3} successfully relaxed the requirement $V_n^2 \to_P V^2$ to the weaker condition $V_n^2 \to_D V^2$. However, when their results are applied to the present martingale array $\{S_{nk}, {\cal F}_{nk}\}_{k\ge 2, n\ge 1}$, a crucial additional assumption imposed in the cited papers, namely,
$
\frac{1}{n}\sum_{j=1}^n H(x_{nj},u_{nj}) \to_P 0,
$
is not satisfied.

The main aim of this paper is to remove this restriction. Under certain regular conditions on $x_{nj}$ and $u_{nj}$, ensuring that $
\frac{1}{n}\sum_{j=1}^n H(x_{nj},u_{nj}) \to_D \int_0^1 \widetilde H(X_t)dt,
$ our main result Theorem \ref {main}  establishes   new asymptotics  involving a functional of stochastic integrals for $S_{nn}$,   providing a framework in an extension to the existing martingale limit theorem. In the remainder of this paper,  our new  result is applied to nonlinear regression with nonstationary time series in Section 3, allowing for the models to include an integrated regressor with certain lags,  which are most commonly used in practice. Using Theorem \ref {main} as a main tool,  Theorem \ref {th8}  investigate the asymptotics
 for the conventional least square estimators in a nonlinear cointegrating regression model.
 These results  improve those in the existing literature.
All technical proofs are given in Section 4.

 Throughout the paper, we denote by $C, C_1, ...$
constants, which may change at each appearance. The notation $\to_D$ ($\to_P$ respectively)
  denotes    convergence in distribution (in probability respectively) for a sequence of random variables (or vectors). The notation $=_D$ denotes equality in distribution.
If $\alpha_n^{(1)}$, $\alpha^{(2)}_n $,..., $\alpha^{(k)}_n$ $(1\le n\le
\infty)$ are random elements on $D[0,1]$ or $D[0, \infty)$ or $\mathbf{R}$,  we will understand the condition
\begin{equation*}
(\alpha_n^{(1)}, \alpha^{(2)}_n,...,\alpha^{(k)}_n)\Rightarrow
(\alpha_{\infty}^{(1)},  \alpha_{\infty}^{(2)},..., \alpha_{\infty}^{(k)})
\end{equation*}
to mean that for all $\alpha_{\infty}^{(1)}$, $\alpha_{\infty}^{(2)}$,...,  $%
\alpha_{\infty}^{(k)}$-continuity sets  $A_1$, $A_2$,...,$A_k,$
\begin{equation*}
P\big(\alpha_n^{(1)}\in A_1, \alpha^{(2)}_n\in A_2,...,\alpha^{(k)}_n\in A_k%
\big)
\to P\big(\alpha_{\infty}^{(1)}\in A_1,  \alpha_{\infty}^{(2)}\in A_2,
...,\alpha_{\infty}^{(k)}\in A_k\big).
\end{equation*}
[see  Theorem 3.1 of \cite{B1}]. As usual,
  $D[0,1]$ ($D[0,\infty)$ respectively) denotes the space of cadlag functions on $[0,1]$ ($[0, \infty)$ respectively), which is equipped with Skorohod topology.

\section{Main results}

We recall that $\{\ep_{nj}, {\cal F}_{nj}\}_{j\ge 1, n\ge 1}$ is a    martingale difference array defined on a common probability space $(\Omega, {\cal F}, P)$ with
  $ \E\big(\ep_{nj}^2 | {\cal F}_{n, j-1} \big)=1,  a.s.$ for all $n, j\ge 1$. For each $n\ge 1$, we always assume that ${\cal F}_{n0}$ is a trivial $\si$-field, $
  {\cal F}_{nj} \subseteq {\cal F}_{n,j+1},  j\ge 1,
 $  and  both random vector
sequences  $\{x_{nj}\}_{j\ge 1}$ (with degree $d_1$) and $\{u_{nj}\}_{j\ge 1}$ (with degree $d_2$) are      adapted to the filtration ${\cal F}_n=\{{\cal F}_{n, j-1}\}_{j\ge 1}$. In particular, we allow for $d_2=\infty$, i.e., we may have   $u_{nj}=(u_{n j, 1}, u_{nj, 2}, ...)$.
This section investigates   the asymptotics of $S_n$ defined by
\bestar
S_n &=&\frac 1{\sqrt n}\sum_{j=1}^n H(x_{nj}, u_{nj})\, \ep_{nj}.
\eestar
To this end,  we make use of the following additional assumptions.

\begin{assump} \la {as2.1}
There exists  a   $d$-dimensional  martingale difference array $\{(\eta_{nj, 1}, ...,\eta_{nj, d}), {\cal F}_{nj}\}\\ _{j\ge 1, n\ge 1}$ on $(\Omega, {\cal F}, P)$ so that $\E\big(\ep_{nj} \eta_{nj, 1}\mid {\cal F}_{n, j-1} \big)=...=\E\big(\ep_{nj} \eta_{nj, d}\mid {\cal F}_{n, j-1} \big)=~0$ for all $n, j\ge 1$ and, on $D_{\mathbb{R}^{1+d+d_1}}[0,\infty)$,
\be
\Big(\frac 1{\sqrt n }\sum_{j=1}^{[nt]} \ep_{nj}, \ \frac 1{\sqrt n}\sum_{j=1}^{[nt]} \eta_{nj},
 \ x_{n, [nt]}\Big) &\Rightarrow& \big(B_t, B_{1t},   X_t\big) , \la {c1.1}
\ee
where $\eta_{nj}=(\eta_{nj, 1}, ...,\eta_{nj, d})$, $ (B_{t}, B_{1t}), t\ge 0,$ is  a  $(1+d)-$dimensional standard Brownian motion and $X=\{X_t\}_{t\ge 0}$ is  a continuous   process adapted to the filtration ${\mathcal F}=\{\si (B_s, B_{1s}, s\le t)\}_{t\ge 0}$.
\end{assump}
\begin{assump} \la {as2.2} $H(x, y)$ is a real function   satisfying the conditions:
\begin{itemize}
\item[(a)] $H(x, y)$ is piecewise continuous  with respect to $x\in \mathbb{R}^{d_1}$;
\item[(b)]  for each  $A>0$,  $\{ g_A^2(u_{nj}) \}_{n\ge 1, j\ge 1}$ is uniformly integrable, where $g_A(y)= \sup_{||x||\le A} |H(x, y)|.$

\end{itemize}
\end{assump}
\begin{assump} \la {as2.3} There exist  piecewise continuous functions $\widetilde H(x)$ and $\widehat H(x)$ on $\mathbb{R}^{d_1}$ and for some sequence $0<m := m_n \to \infty$ with ${n}/{m} \to\infty$ so that, for each fixed $x$,
\be
\max_{m\le k\le n-m} \E \,\Big|\frac 1m \sum_{j=km+1}^{(k+1)m}\, H(x, u_{n j})-\widetilde H(x)\,\Big|
&\to& 0,  \la {c1.2} \\
\max_{m\le k\le n-m} \E\, \Big|\frac 1m \sum_{j=km+1}^{(k+1)m}\, H^2(x, u_{n j})-\widehat H^2(x)\,\Big|
&\to& 0,  \la {c1.2a}
\ee
as $m\to \infty$.
\end{assump}

We have the following  main result.

\begin{thm} \la {main} Suppose that Assumptions  \ref {as2.1} - \ref {as2.3} hold. Then,  as $n\to \infty$, we have
\be
  S_n&\to_D &S:=\,\int_0^1 \widetilde H(X_t)dB_{t}+\si(X)\, {\cal N} , \la {m1}
\ee
where $ \si^2(X)=\int_0^1 \big[\widehat H^2(X_t)- \widetilde H^2(X_t)\big] dt  $  and ${\cal N}$ is a standard normal variate independent of $B=\{B_t\}_{t\ge 0}$ and  $X=\{X_t\}_{t\ge 0}$.
\end{thm}

\begin{rem} \la {r2.1} The condition (\ref {c1.1}) in Assumption  \ref {as2.1} is standard in investigating general limit theorems for martingale arrays. See, for instance, \cite{W1, W2}. The    martingale difference array $\{\eta_{nj}, {\cal F}_{nj}\}_{j\ge 1, n\ge 1}$ is usually involved in exploring the limit law of   $\{x_{nj}\}_{j\ge 1}$.  Its existence is obvious in many applications as seen in the following (\ref {c1.1-01}).  If  $X_t$ is measurable with respect to $\sigma(B_s: s \le t)$, we may take $\eta_{nj} \equiv 0$.  Note that,  if $X$ is  independent of $B$, we may rewrite (\ref {m1}) as
$
S_n\to_D \,\si_1(X)\, {\cal N},
$
where $\si_1^2(X)=\int_0^1 \widehat H^2(X_t) dt $.
Assumption  \ref {as2.2} is weak, where  the piecewise continuity of $H(x, y)$ in Assumption  \ref {as2.2}(a) is enough for many empirical applications. Assumption  \ref {as2.2}(b) suggests that   (\ref {c1.2}) and (\ref {c1.2a})  only involve the  stationary components of $S_n$, which are easy to be verified in general. In particular, when $u_{nj}\equiv u_j$  is (strictly) stationary  satisfying  $E g^2_A(u_1) <\infty$ for each $A>0$, the  conditions (\ref {c1.2}) and (\ref {c1.2a}) hold with $ \widetilde H(x)\, =\E H(x, u_{1})$ and $ \widehat H^2(x)\, =\E H^2(x, u_{1})$, respectively.
\end{rem}

\begin{rem} \la {r2.2a} For any square integrable martingale difference array $\{X_{nj}, {\cal F}_{nj}\}_{j\ge 1, n\ge 1}$, we may write
$
X_{nj} =\si_{nj}\, \ep_{nj}$, where $\si_{nj}^2=\E\big(X_{nj}^2 | {\cal F}_{n, j-1} \big) $ and
$\{\ep_{nj}, {\cal F}_{nj}\}_{j\ge 1, n\ge 1}$ is a    martingale difference array  with
  $ \E\big(\ep_{nj}^2 | {\cal F}_{n, j-1} \big)=1,  a.s.$ for all $n, j\ge 1$. Theorem \ref {main} imposes the structure $\si_{nj}= H(x_{nj}, u_{nj})$ with that $x_{nj}$ is  a sequence of  random  vectors satisfying $x_{n, [nt]}\Rightarrow X_t$ on $D_{\mathbb{R}^{d_1}}[0,\infty)$.  In earlier works, a similar structure $\si_{nj}= H(r_n \,x_{nj}, u_{nj})$, where $r_n\to \infty$ is a sequence of constants,  was considered in \cite{W1}.  We mention that  the asymptotics of $\sum_{j=1}^nH( x_{nj}, u_{nj})\ep_{nj}$ are essentially different from those of   $\sum_{j=1}^nH( r_n x_{nj}, u_{nj})\ep_{nj}$ with $r_n\to \infty$.  In latter case, the limiting distribution is usually a mixture of normal distributions rather than a function of stochastic integrals. Indeed, under certain regular conditions, when $r_n/n\to 0$ and $r_n\to \infty$, we may have
  \bestar
  \big(\frac {r_n}n\big)^{1/2}\sum_{j=1}^nH( r_n x_{nj}, u_{nj})\ep_{nj}\to_D c_0\,L_X^{1/2}(1, 0) {\cal N}
  \eestar
  where $c_0$ is a constant, $L_X^{1/2}(1, 0)$ is a local time process of $X$ and ${\cal N}$ is a normal variate independent of $X$. For more details in related topics, together with the applications in financial econometrics, we refer to \cite{W2}, \cite{WP1, WP2, WP3}, \cite{WPK1} and references therein.

\end{rem}
\begin{rem}\la {r2.2}
The limiting result (\ref {m1}) of $S_n$ depends essentially on its conditional variance $V_n^2:=V_{nn}^2,$ where, for $k\ge 1$,
\bestar
V_{nk}^2 &=& \frac 1n  \sum_{j=1}^{k} H^2(x_{nj}, u_{nj})
\eestar
and the conditional co-variance $R_{nn}$ between $S_n$ and $\frac 1{\sqrt n }\sum_{j=1}^{[nt]} \ep_{nj}$, where, for $k\ge 1$,
\bestar
R_{nk} &=& \frac 1n  \sum_{j=1}^{k} H(x_{nj}, u_{nj}) .
\eestar
Under   (\ref {c1.2}) and (\ref {c1.2a}),   Lemma \ref {lem3} in Section \ref {sec4} yields the results:
\bestar
 \max_{1\le k\le n} \big| R_{nk}-\frac 1n \sum_{j=1}^{k} \widetilde H(x_{nj})\big| &=&o_P(1), \\
  \max_{1\le k\le n} \big| V_{nk}-\frac 1n\sum_{j=1}^{k} \widehat  H^2(x_{nj})\big| &=&o_P(1) ,
\eestar
indicating that
$
R_{nn} \to_D \int_0^1   \widetilde  H(X_t) dt $ and $V_{n}^2 \to_D \int_0^1  \widehat  H^2(X_t) dt.$
These facts play a central role in the proof of Theorem \ref {main}, where we  have in fact established the following joint convergence: on $D_{\mathbb{R}^{d+d_1+3}}[0, \infty)$,
\be
&& \Big(\frac 1{\sqrt n }\sum_{j=1}^{[nt]} \ep_{nj}, \ \frac 1{\sqrt n}\sum_{j=1}^{[nt]} \eta_{nj},
 \ x_{n, [nt]}, V_{n}^2, S_n\Big) \no\\
&&\Rightarrow\ \big(B_t, B_{1t},   X_t, \int_0^1  \widehat  H^2(X_t) dt,  S\big). \la {ad35}
\ee
\end{rem}

\medskip
We next consider  a simple corollary  of Theorem \ref {main}, illustrating the difference from previous works. More general application will be investigated in Section \ref {sec3}. 

Let $\ep_i$ be iid random variables with $E\ep_1=0$ and $E\ep_1^2=1$. Let $x_{nj}=\frac 1{\sqrt n} \sum_{k=1}^{j-1}\ep_k$ and write
\bestar
S_{n1} &=&\frac 1{\sqrt n} \sum_{j=1}^n x_{nj} K(\ep_{j-1})\ep_j, \quad
 S_{n2} \ =\ \frac 1{\sqrt n} \sum_{j=1}^n\big[ x_{nj}+ K(\ep_{j-1})\big]\, \ep_j, 
\eestar
where $K(y)$ is a real function  satisfying that $|K(y)|\le C\, |y|^{\al}$ for some $\al\ge 0$.  By Theorem \ref {main} with $u_{nj}=K(\ep_{j-1})$, $=1$ and $H(x, y)=xK(y)$ or $H(x, y)=x+K(y) $ respectively, we have the following corollary.
\begin{cor} Suppose that $\E |\ep_1|^{\max\{2\al, 2\}}<\infty$. Then, as $n\to\infty$, we have
\be
S_{n1} &\to_D&  \E K(\ep_1)\, \int_0^1 B_tdB_t + \si_K\, \sqrt  {\int_0^1 B_t^2 dt}\, {\cal N}, \\
S_{n2} &\to_D&\int_0^1 \big[B_t+ \E K(\ep_1)\big] dB_t + \si_K\, {\cal N},
\ee
where $\si_K^2= var \big[K(\ep_1)\big]$, $B_t$ is a standard Brownian motion and ${\cal N}$ is a standard normal variate independent of $B=\{B_t\}_{t\ge 0}$.
\end{cor}

If $\E K(\ep_1)=0$, the limiting distribution of $S_{n1}$ is mixture normal, which is similar to that provided in recent works such as \cite{J1},  \cite{W1, W2} and  \cite{RV1, RV2, RV3}. When $\E K(\ep_1)\not=0$, the limiting distributions of $S_{n1}$ and $S_{n2}$ involve a functional of stochastic integrals, which seems new to the literature.

\section{A regression application } \la {sec3}

This section considers nonlinear cointegrating regression model having the following form:
\be
y_t=f(x_t, \Delta x_{t}, ..., \Delta x_{t-d}, \  \theta)+u_t, \quad t=d, d+1, ... , \la {m2g}
\ee
for some $d\ge 0$, where $f(. ; ..., \theta)$ is a given real function indexed by $\theta=(\theta_1,..., \theta_q) $, a vector of unknown parameters, $x_t$ is an integrated   regressor, $\Delta x_k=x_{k}-x_{k-1}, k=t, ..., t-d$ denotes  the lags of the integrated regressor $x_t$, and $(u_t, {\cal F}_t)_{t\ge 1}$ is a sequence of martingale differences such that $x_t$
is adapted to ${{\cal F}_{t-1}}$. Let $\Theta$ be a compact set of $R^q$ and assume the unknown parameters $\theta\in \Theta$. The least squares  estimator (LSE) $\hat \theta_n$ of $\theta$, an interior of $\Theta$,  is defined by
\bestar
\widehat \theta_n =\arg\min_{\theta\in \Theta}\, \sum_{t=1}^n \big[y_t-f(x_t, \Delta x_{t}, ..., \Delta x_{t-d},\ \theta)\big]^2.
\eestar The asymptotics of $\hat \theta_n$ without  lags $\Delta x_{t}, ..., \Delta x_{t-d}$ were initially considered in \cite{PP1} and then later by  \cite{CPP1},  \cite{dH1}, and \cite{CW1}.
More recently, \cite{W3} investigated  the asymptotics by allowing for a wide class of commonly used volatility models such as GARCH, time-varying GARCH, and nonlinear GARCH. See also \cite{JW1} for the  related works on  weighted  nonlinear cointegrating regression.
Using Theorem \ref {main}, this section will explore the asymptotics of $\hat \theta_n$ by allowing for
  certain lags $\Delta x_{t}, ..., \Delta x_{t-d}$   in  the  model  (\ref {m2g}).


The remainder of this section is organized as follows. Section 3.1 presents the assumptions on $x_t$ and $u_t$, which are similar to those in \cite{W3} and encompass widely used volatility models. Section 3.2 presents the asymptotic results.

\subsection{Assumptions on $x_t$ and $u_t$}

Throughout the section, let
 $\lambda_{i}\equiv (\ep_{i-1},\nu_i), i\in \mathbb{Z}$, be a sequence of i.i.d. random vectors
with ${E} \lambda _{0}=0$, $E\ep_0^2=E\nu_0^2=1$ and $\rho=E\ep_0\nu_1$, and $\{\lambda_{i}^*\}_{i\in \mathbb{Z}}$ be an independent copy of $\{\lambda_{i}\}_{i\in \mathbb{Z}}$.
We make use of the following assumptions on $x_t$ and $u_t$ in the asymptotic development.

\medskip
\begin{assump} \la {as3.1}
 $x_t= x_{t-1}+\xi_t,$
 where  $ \xi_t=\sum_{j=0}^{\infty}\phi_j\nu_{t-j}$. The coefficients $\phi _{j},j\geq 0,$ satisfy one of the following conditions:
\begin{quote}
\begin{itemize}
\item[\quad\textbf{LM.}] $\phi _{j}\sim j^{-\mu }\,l(j), 1/2<\mu <1$ and $
l(k)$ is a function slowly varying at $\infty $.
\item[\quad \textbf{SM.}] $\sum_{j=0}^{\infty }|\phi _{j}|<\infty $ and $\phi \equiv
\sum_{j=0}^{\infty }\phi_{j}\not =0$.
\end{itemize}
\end{quote}
\end{assump}

\begin{assump} \la {as3.2}  $u_t=\si_t\ep_{t}$ with   $\si_t^2=\si(\lam_t, \lam_{t-1},...)$, where $\si(...)$ is a measurable function satisfying $ \E\si^2( \lam_0, \lam_{1},...)<\infty$ and
    \be
 E\big|\si( \lam_m, \lam_{m-1},...) -   \si( \lam_m, \lam_{m-1},...,\lam_1, \lam_0^*, \lam_1^*,...)\big|^2\ \le\ C \,m^{-\al},\la {er1a}
\ee
for any $m\ge 1$ and some $\al >\left \{\begin{array}{ll}
4/(2\mu-1), & \mbox{under {\bf LM}}, \\
4, & \mbox{under {\bf SM}.}
\end{array}
\right.
$
\end{assump}

Assumption \ref {as3.1} allows for short memory (under {\bf SM}) and long memory (under {\bf LM}) innovations driving the  integrated regressor $x_k=\sum_{j=1}^k\xi_j$, which is quite general in practice. Define
\be
d_{n}^{2} &=&\mathbb{E}|\sum_{k=1}^{n}\xi_{k}|^{2}\ \sim\ \left \{\begin{array}{ll}
c_{\mu }\,n^{3-2\mu }l^{2}(n), & \mbox{under {\it LM}}, \\
\phi^2\, n, & \mbox{under {\it SM},}
\end{array}%
\right. \la {f0}
\ee
where $c_{\mu}$ is a constant.
 Standard functional limit theory (see \cite{BC1} or Theorem 2.21 in \cite{W2} with a minor modification) shows that
\be
&& \big(\frac{1}{\sqrt{n}}\sum_{i=1}^{\lfloor nt\rfloor }\epsilon _{i},\ \frac{1
}{\sqrt{n}}\sum_{i=1}^{\lfloor nt\rfloor }\nu _{i},\ \frac{1}{\sqrt{n}}
\sum_{i=1}^{\lfloor nt\rfloor }\nu _{-i}, \frac {1}{d_n}x_{[nt]}\big) \no\\
&&\Rightarrow\ \big(U_{1t},\
U_{t},\ U_{-t},\  X_t\big), \la {f2}
\ee
on $D_{\mathbb{R}^{4}}[0,\infty)$, where $(U_{1t},U_{t})_{t\ge 0}$ is a bivariate Brownian motion with covariance matrix:
$
\Omega =
\begin{pmatrix}
1 & \rho \\
\rho & 1
\end{pmatrix}
,
$
 $\{U_{-t}\}_{t\ge 0}$ is an independent copy of $\{U_t\}_{t\ge 0}$ and
 $X_t$ is defined by
\bestar
X_t &=&\left \{\begin{array} {ll}
U_{3/2-\mu}(t), & \mbox{under {\bf LM}},\\
U_{1/2}(t), & \mbox{under {\bf SM}},
\end{array}\right.
\eestar
where $U_H(t)$ is a fractional Brownian motion  that has the following  representation: with $a_+=\max\{a, 0\}$,
\bestar
U_{H}(t) &=& \frac {1}{\Gamma(H+1/2)}\int_{-\infty}^t(t-s)_{+}^{H-1/2}-(-s)_+^{H-1/2} dU_s.
\eestar

In terms of (\ref {f2}), there exists   a    martingale vector difference array $\{(\eta_{1j}, \eta_{2j}), {\cal F}_{j}\}_{j\ge 1}$, where ${\cal F}_{j}=\si(\lam_i, i\le j)$ so that $\E\big(\ep_{j} \eta_{1j}\mid {\cal F}_{ j-1} \big)=\E\big(\ep_{j} \eta_{2j}\mid {\cal F}_{ j-1} \big)~0$ for all $ j\ge 1$ and, on $D_{\mathbb{R}^{4}}[0,\infty)$,
\be
\Big(\frac 1{\sqrt n }\sum_{j=1}^{[nt]} \ep_{j}, \ \frac 1{\sqrt n}\sum_{j=1}^{[nt]} \eta_{j},
 \ \frac {1}{d_n}x_{[nt]}\Big) &\Rightarrow& \big(B_t, B_{1t},   X_t\big) , \la {c1.1-01}
\ee
where  $ (B_{t}, B_{1t}), t\ge 0,$ is  a  $3-$dimensional standard Brownian motion and $X=\{X_t\}_{t\ge 0}$ is  a continuous   process adapted to the filtration ${\mathcal F}=\{\si (B_s, B_{1s}, s\le t)\}_{t\ge 0}$. Indeed, if  we  take  $\eta_{2j}= \nu_{-j}$ and
$
\eta_{1j} =\big(\rho\ep_j- \nu_j\big)/(1-\rho),
$
it is routine to see that
\bestar
\E\big(\ep_{j} \eta_{1j}\mid {\cal F}_{ j-1} \big)=\E\big(\ep_{j} \eta_{2j}\mid {\cal F}_{ j-1} \big)=\E\big(\ep_{j} \eta_{1j}\big)=0
\eestar and (\ref {c1.1-01}) is true by (\ref {f2}) and the continuous mapping theorem.

Under the settings of Assumption \ref {as3.2}, we have that $(u_k, {\cal F}_{k})_{k\ge 1}$, where ${\cal F}_k$ is an $\si$-field generated by $\lam_{k+1}, \lam_{k},...,$ forms a martingale difference with $E(u_k^2|{\cal F}_{k-1})=\si_k^2=\si(\lam_k, \lam_{k-1},...)$. The error process $u_k$ having a martingale difference structure that is similar to that used in \cite{W3}, indicating that (\ref {m2g}) allows for a wide class of models of  heteroscedasticity such as   GARCH  and nonlinear GARCH models. The following exam from \cite{W3} provides an illustration. More related  results can be found in  \cite{W3},  \cite{WM1}, \cite{SW1}, and \cite{PW1}.

\begin{exap} (Nonlinear GARCH model) Let $u_t=\si_t^{1/2} \ep_t$, where $\si_t$ is defined by
\bestar
\si_t= \al_0+ \al\, F(\si_{t-1})+\beta \si_{t-1}+  \ga u_{t-1}^2.
\eestar
This is a nonlinear GARCH model introduced by \cite{LS1}. Suppose that a constant $A_0$ exists such that
$|F(x)-F(x')|\le A_0|x-x'|$ for all $x, x'\in R$. Since we may rewrite
$
\si_t=\al_0 +R(\si_{t-1}, \ep_{t-1}),
$ where $R(y, \ep)=\al F(y)+\beta y+\ga \ep^2y$, it follows from Proposition 3.1 of \cite{W3} that
 $u_t$ satisfies Assumption \ref {as3.2} if $A_0>0, \al\ge 0, \beta\ge 0, \ga\ge 0$ satisfying $\big[\E\big(\al A_0+\beta+\ga \ep_0^2\big)^2\big]^{1/2}<1$.

\end{exap}

\subsection{Asymptotics of $\hat \theta_n$}
This section establishes the limiting distribution of $\hat{\theta}_n$ under certain regular conditions on $f(x, y, \theta)$. To this end,  define $\dot{f}(x, y, \theta)= (\dot{f}_1,...,\dot{f}_q)^{\prime},$ where  $\dot{f}_i=\frac {\partial f(x, y, \theta)}{\partial \theta_i}, i=1,..., q$, and write $p(x,  y, \theta)$ for either  $f(x, y, \theta)$ or one of  $\dot{f}_i, i=1,..., q$. We need the following restrictions  on $p(x, y, \theta)$.

\begin{assump}\la {as3.3} There exist a real continuous function $T_p(x)$ and a constant $\beta \geq 0$ such that
\begin{enumerate}
	\item[(i)] for each $\theta_1,\theta_2 \in \Theta$, $(x, y)\in {\mathbb R}^{1+d}$ and for some $0<\al\le 1$,
\be
|p(x, y, \theta_1) - p(x, y, \theta_2)| &\le& ||\theta_1 - \theta_2||^{\al}\,T_p(x)(1+||y||^\beta); \la {condy}
 \ee

	\item[(ii)] for any bounded $x$,
\be
p(l x, y, \theta) &=& v_p(l)\, h_p(x, y, \theta)\, +\, R_p(l x, y, \theta), \la {3.14}
\ee
where $v_p(l)$ is a positive real function that is bounded away from zero as $l$ is large, $ h_p(x,y,  \theta)$ for each $\theta\in \Theta$ is a continuous function with respect to $x$ satisfying $ |h_p(x,y,  \theta)|\le C_x (1+||y||^\beta)$  and
\bestar
\sup_{\theta\in \Theta, y\in R}\frac {|R_p(lx, y,  \theta)|}{T_p(l x)(1+||y||^\beta)} =o(1),\quad \mbox{as $l\to \infty$.}
\eestar

\item[(iii)] $T_p(l x)\le v_p(l)\,(1+|x|^\ga) $ for some $\ga>0$,  as $|l x|\to \infty $.

 \end{enumerate}

\end{assump}

To introduce the main result in this section, let $\theta_0$ be   an interior of $\Theta$ and we make use of the following notations:
 \bestar
\widetilde h_f (x, \theta) &=& \E h_f(x, \xi_0, ..., \xi_{-d}; \theta),\qquad  \widetilde \si_f (x, \theta) \ =\ \E \big[h_f(x, \xi_0, ..., \xi_{-d}; \theta) \si_0\big],\\
  { H_{\dot f}} (x, \theta)&=&\big( h_{\dot{f}_1}(x,\xi_0\,...\, \xi_{-d};\theta),\,..., \,h_{\dot{f}_q}(x,\xi_0, ..., \xi_{-d};\theta)\big),\\
\widetilde  { h_{\dot f}} (x, \theta)&=&\E  { H_{\dot f}} (x, \theta), \qquad \widetilde  { \si_{\dot f}} (x, \theta)=\E \big[ { H_{\dot f}} (x, \theta) { H_{\dot f}} (x, \theta) ^{\prime}\si_0^2\big].
\eestar
It is routine from Assumption \ref {as3.3} that these notations are well defined if
 \be
 \E \big(|\xi_{-j}|^{2\beta}\si_0^2\big)<\infty, \quad \mbox{ for $j=0, -1, ..., -d$.} \la {addcon}
\ee

\ignore{
\begin{assump}\la {ass3.4}
For each $\eta>0$, we have
\begin{itemize}
\item[(a)]
$
\int_{|x|\le \eta}\big[\widetilde h_f(x, \theta)-\widetilde h_f(x, \theta_0)\big]^2dx \not=0,\ \mbox{for any $\theta\not=\theta_0$}; and
$
\item[(b)]
$
\int_{|x|\le \eta}\widetilde  { h_{\dot f}}(x, \theta_0)\,  \widetilde  { h_{\dot f}}^{\prime}(x, \theta_0)dx \ \ \mbox{is a positive-definite matrix}. 
$
\end{itemize}
\end{assump}
}

\begin{thm} \la{th8} In addition to Assumptions \ref {as3.1}-\ref{as3.3}, suppose that  (\ref {addcon}) holds and
\be
&& \int_{|x|\le \eta}\big[\widetilde h_f(x, \theta)-\widetilde h_f(x, \theta_0)\big]^2dx \not=0,\ \mbox{for any $\theta\not=\theta_0$,}\hskip 5cm \la {tc1} \\
&& \int_{|x|\le \eta}\widetilde  { h_{\dot f}}(x, \theta_0)\,  \widetilde  { h_{\dot f}}^{\prime}(x, \theta_0)\, dx \ \ \mbox{is a positive-definite matrix}, \la {h12}
\ee
for each $\eta>0$.
 Then, as $n\to\infty$,
\be \la{thmHomoLimit1.eqn1}
D_n\, ( \hat{\theta}_n - \theta_0) &\rightarrow_D&\,\Big( \int_0^1 \Psi(t) \Psi(t)' dt \Big )^{-1} \,\Big[
\int_0^1 \widetilde G_t\,dB_t+ {\Omega}^{1/2}\, {\cal N}_q\Big],
\ee
where  $D_n=\sqrt n\, \mbox{diag}\, \big(v_{\dot {f}_1}(d_n), ...,v_{\dot {f}_q}(d_n)\big)$,
\bestar
\Psi(t)  &=& \widetilde  { h_{\dot f}}(X_t, \theta_0), \quad \widetilde G_t= \widetilde \si_f(X_t, \theta_0), \quad \widehat G_t =\widetilde  { \si_{\dot f}}  (X_t, \theta_0) \no\\
\Omega &=&\int_0^1\big[\widehat G_t-\widetilde G_t\widetilde G_t^{\prime}\big]dt
\eestar and   ${\cal N}_q$ is an $q$-dimensional standard normal vector independent of $B$ and  $X
$.

\end{thm}

\begin{rem} As indicated in \cite{CW1} and \cite{PP1}, nonlinear cointegrating regressions with structures described in Theorem \ref{th8} are useful for modeling money demand functions. In such cases, $y_t$ is the logarithm of the real money balance, $x_t$ is the nominal interest rate, and $f$ can either be $f(x, \al, \beta) = \al + \beta \log | x |$ or $f(x, \al, \beta) = \al + \beta \log (\frac{1 + |x|}{|x|}).$ See \cite{Bd1} and \cite{BKO2} for empirical investigations of the estimation of money demand functions in the United States and Japan, respectively. See also \cite{BKO1} for the derivation of these functional forms from the underlying money demand theories studied in macroeconomics.
\end{rem}

\section{Proofs} \la {sec4}


Except mentioned explicitly, we use the same notations as in previous sections.
We start with the proof of Theorem \ref {main}.  Three lemmas that support the proof of
Theorem \ref {main} are given in Subsections \ref {subs1} -- \ref {subs3}, respectively.  The proof   of Theorem \ref {th8} is given in Section \ref {sec4.2}.
Without loss of  generality, in the proof of Theorem \ref {main}, we assume that $d_1=1$ and $H(x, y)$ is continuous with respect to $x$. General extension to $d_1>1$ and the piecewise continuity of $H(x,y)$ is obvious from the proof with some notation changes, and hence the details are omitted.

\subsection{Proof of Theorem \ref {main}}
We first prove (\ref {m1}) under an additional restriction:
\begin{itemize}
\item[{\bf R}:]$H(x, y)$ is continuous  with respect to $x$ and $H(x, y)=0$ for all  $|x|\ge N$, where $N\ge 1$ is a constant.  \la {addic}
\end{itemize} This  restriction will be removed later. Note that, under this additional restriction, the continuous functions $\widetilde H(x)$ and $\widehat H(x)$ given in Assumption \ref {as2.3} can be chosen so that  $\widetilde H(x)=\widehat H(x)=0$ for all $|x|\ge N.$  This fact will be used in the proof and the following three lemmas without further explanation.
Write
\bestar
S_{1n} &=& \frac 1{\sqrt n} \sum_{j=1}^n \widetilde H (x_{nj})\,  \ep_{nj} \quad \mbox{and}\quad S_{2n} = \frac 1{\sqrt n} \sum_{j=1}^n K_j(x_{nj}, u_{nj}) \, \ep_{nj},
\eestar
where $K_j(x, y)=H(x, y)-\widetilde H (x)$. It suffices to show that
\be
\big( S_{1n}, \ S_{2n}\big) &\to_D& \Big(\int_{0}^1  \widetilde H(X_t)dB_t, \ \si(X)\, {\cal N } \Big),  \la {3.1}
\ee
where ${\cal N } $ is a standard normal variate independent of $(B_t, X_t)_{t\ge 0}$ or, equivalently, independent of $(B_t, B_{1t})_{t\ge 0}$ since $X=\{X_t\}_{t\ge 0}$ is  adapted to the filtration ${\mathrm F}=\{\si (B_s, B_{1s}, s\le t)\}_{t\ge 0}$.
To this end, for $j, k=1, 2, 3,..., n$,  let $y_{nj} = \frac 1{\sqrt n}  K_j(x_{nj}, u_{nj}) \, \ep_{nj}, $
\bestar   V_{nk}^2=\frac 1n\sum_{j=1}^k  K_j^2(x_{nj}, u_{nj})  ,\quad  \widetilde V_{nk}^2=\frac 1n\sum_{j=1}^k  \big(\widehat H^2(x_{nj})- \big[\widetilde H(x_{nj})\big]^2\big).
\eestar
Suppose that $l_j, j\ge 1,$ are iid $N(0,1)$ random variables and $\{l_j\}_{j\ge 1}$ is  independent of all other variables $\ep_{nj}, u_{nj}, $ and $x_{nj}$. We further write
\bestar
\widehat y_{nj} = \left\{\begin{array}{ll}
y_{nj}, &\mbox{if $j\le n$}, \\
\frac 1{\sqrt n}l_i, & \mbox{if $j>n$,}
\end{array}\right. \quad \widehat V_{nk}^2 &=& \left\{\begin{array}{ll}
V_{nk}^2, &\mbox{if $k\le n$}, \\
V_{nn}^2+\frac 1n \sum_{j=n+1}^k l_j^2, & \mbox{if $k>n$,}
\end{array}\right.
\eestar
and   $\tau_n(t)=\inf\{k:  \widehat V_{nk}^2\ge t\}$ for $t> 0$. We mention that such a  definition of $\widehat V_{nk}^2$ ensures that  $ \widehat V_{nk}^2\to\infty, a.s.$ for each fixed $n$. We introduce the following  three lemmas. Their proofs will be given in  Subsections \ref {subs1} -- \ref {subs3}, respectively.

\begin{lem} \la {lem31}
If,  in addition to  Assumptions \ref {as2.1} - \ref {as2.3}, $H(x, y)=0$ for all  $|x|\ge N$ where $N\ge 1$ is a constant, then
\be
\max_{1\le j\le n} |y_{nj}| &=&o_P(1).\la {4.6a}
\ee

\end{lem}

\begin{lem} \la {lem3} If,  in addition to   Assumptions \ref {as2.1} -  \ref {as2.3}, $H(x, y)=0$ for all  $|x|\ge N$ where $N\ge 1$ is a constant, then
\be
&& \frac 1 {\sqrt n} \max_{1\le k\le n} \big|\big| \sum_{j=1}^{k} \E \big(\xi_{nj} y_{nj} |{\cal F}_{n, j-1}\big)\big|\big|  \no\\
&= &\frac 1n \max_{1\le k\le n} \big| \sum_{j=1}^{k}\big[H(x_{nj}, u_{nj}) -\widetilde H(x_{nj})\big]\big| = o_P(1), \la {4.19}
\ee
where $\xi_{nj}=\big(\,\ep_{nj}, \,\eta_{nj}\big)$, and
\be
\max_{1\le k\le n}\big|V_{nk}^2-\widetilde V_{nk}^2\big|&=& o_P(1).\la {3.4}
\ee
\end{lem}

\begin{lem} \la {lem4} If,  in addition to  Assumptions \ref {as2.1} -  \ref {as2.3}, $H(x, y)=0$ for all  $|x|\ge N$ where $N\ge 1$ is a constant, then
\be
\sum_{j=1}^{\tau_n(t)} \widehat y_{nj}& \Rightarrow& B'_t \quad \mbox{on $D_{R}[0, \infty)$,} \la {1.12y}
\ee  where  $B'=\{B_t'\}_{t\ge 0}$ is a standard Brownian motion. Furthermore, on $D_{\mathbb{R}^{4}}[0, \infty)$, we have
 \be
Z_n(t)&:=& \Big(\frac 1{\sqrt n }\sum_{j=1}^{[nt]} \ep_{nj}, \ \frac 1{\sqrt n}\sum_{j=1}^{[nt]} \eta_{nj},
 \ x_{n, [nt]},\ \sum_{j=1}^{\tau_n(t)}\widehat  y_{nj}\Big) \no\\
 &\Rightarrow& \big(B_t,\ B_{1t}, \  X_t,\ B_t'\big)  ,
\la {1.12a}
\ee
where $B'$ is independent of $X=\{X_t\}_{t\ge 0}$ and ${\mathrm F}=\{\si (B_s, B_{1s}, s\le t)\}_{t\ge 0}$.
\end{lem}

We now turn back to the proof of (\ref {3.1}).
In fact, in terms of (\ref {1.12a}),   the same arguments as in \cite{KP1} (also, see,  \cite{JS1}) yield that
\be
&& \Big(\frac 1{\sqrt n }\sum_{j=1}^{[nt]} \ep_{nj},\ \frac 1{\sqrt n}\sum_{j=1}^{[nt]} \eta_{nj}, \   x_{n, [nt]}, \  S_{1n}, \ \widetilde V_{nn}^2, \ \frac 1{\sqrt n} \sum_{j=1}^{\tau_n(t)}\widehat  y_{nj}\Big) \no\\
&\Rightarrow& \(B_{t}, \ B_{1t}, \ X_t, \ \int_{0}^1 \widetilde H(X_{t})dB_{t},\ \si^2(X), \ B_t' \).
\la {3.7}\ee
Write $\tau_n^{-1}(k)=\inf\{t\ge 0: \tau_n(t)\ge k\}$. Note that \bestar
\widehat H^2(x_{nn})+ \big[\widetilde H(x_{nn})\big]^2\to_D\widehat H^2(X_1)+ \big[\widetilde H(X_1)\big]^2 =O_P(1),
\eestar
by using $x_{n, [nt]}\Rightarrow X_t$ given Assumption \ref {as2.1},
It follows from  (\ref {3.4}) that
\be
|\tau_n^{-1} (n)- \widetilde V_{nn}^2| &\le&
|\tau_n^{-1} (n)- \widehat V_{nn}^2| + |V_{nn}^2-\widetilde V_{nn}^2|\no\\
&\le&   V_{nn}^2- V_{n, n-1}^2+\frac 1n l_{n+1}^2 + |V_{nn}^2-\widetilde V_{nn}^2|\no\\
&\le & 3\max_{1\le k\le n}\big|V_{nk}^2-\widetilde V_{nk}^2\big| + \frac 1n \big(\widehat H^2(x_{nn})- \big[\widetilde H(x_{nn})\big]^2\big)+\frac 1n l_{n+1}^2\no\\
&\to_\mathrm{P}&0. \la {a2.3}
\ee
 This, together with  (\ref {3.7}), yields that
 \be
 \Big( S_{1n},  \ \tau_n^{-1} (n), \ \ W_n(t)\Big)
&\Rightarrow&\(\int_{0}^1  \widetilde H(X_{t})dB_{t},\ \si^2(X),\ B_{t}^\prime\), \la {kk1}
\ee
on $D_{\mathrm{R}^{3}}[0,\infty)$, where $W_n(t)=\frac 1{\sqrt n}\sum_{j=1}^{\tau_n(t)} \widehat y_{nj}.$ Let $D_0$ be the subset of all functions $V=\{V(t)\}_{t\ge 0}$ in $ D_{\mathrm{R}}[0, \infty)$ that are nondecreasing with $V(0)\ge 0$.
Since $W_n(t)\in D_{R}[0, \infty)$ and  the composition map $\psi (U, V)$ from $ D_{\mathrm{R}}[0, \infty)\times D_0$ to $D_{\mathrm{R}}[0, \infty)$ defined by
\bestar
\psi (U, V): =U(V(t)), \quad t\ge 0,
\eestar is continuous at  $(B', \si(X))$ (see, \cite{B1}, for instance), it  follows from  (\ref {kk1}) and the continuous mapping theorem  that
\bestar
&&( S_{1n},\ \tau_n^{-1} (n) , \ W_n\big[\tau_n^{-1} (n)\big] \Big)\no\\
&\to_D&\(\int_{0}^1  \widetilde H(B_{1t})dB_{t},\ \si^2(X), \ B_{\si^2(X)}^\prime\)=_D \(\int_{0}^1 \widetilde H(B_{1t})dB_{t},\ \si^2(X),\ \si(X)\, {\cal N}\),
\eestar
where  ${\cal N}$ is a standard normal variate independent of $\{B_t, X_t\}_{t\ge 0}$ and we have used the fact that
 $X=\{X_t\}_{t\ge 0}$ is  adapted to  $F=\{\si (B_s, B_{1s}), s\le t\}_{t\ge 0}$  and $B^\prime$ is independent of $F$. Now, by noting
\bestar
\big|S_{2n}-  \frac 1{\sqrt n} \sum_{j=1}^{\tau_n[\tau_n^{-1} (n)]} \widehat y_{nj}\big|\le  \big|y_{nn}\big|\to_\mathrm{P}0
\eestar
owe to (\ref {4.6a}), we have
\bestar
&& \Big( S_{1n}, \ \tau_n^{-1} (n), \  S_{2n}\Big) =( S_{1n},\ \tau_n^{-1} (n) , \ W_n\big[\tau_n^{-1} (n)\big] \Big)+o_P(1)\no\\
&\to_D& \(\int_{0}^1 \widetilde H(B_{1t})dB_{t},\ \si^2(X),\ \si(X)\, {\cal N}\),
\eestar
where  ${\cal N}$ is a standard normal variate independent of $\{B_t, X_t\}_{t\ge 0}$.
This proves (\ref {3.1}) and hence (\ref {m1}) under the additional restriction {\bf R.}

To remove this restriction {\bf R},  for any $K>0$, we assume that  $H_K(x, y)=H(x, y)\xi_K(x)$,   $\widetilde H_K(x)=\widetilde H (x)\xi_K(x)$ and $\widehat H_K(x)=\widehat H (x)\xi_K(x)$
with
\bestar
\xi_K(x)=\left\{
\begin{array}{ll}
1, & |x|\le K,\\
2-|x|/K, & K<|x|<2K,\\
0, & |x|\ge 2K.
\end{array}
\right.
\eestar
It is readily seen that,  for each $K>0$,  $H_K(x, y)$ is a continuous function with respect to $x$,  satisfying
 Assumption \ref {as2.2}, $H(x, y)=0$ for all  $|x|\ge N=2K$,  and Assumption \ref {as2.3} with  $\widetilde H(x)$ and $\widehat H(x)$ being replaced by $\widetilde H_K(x)$ and $\widehat H_K(x)$, respectively. Therefore, it follows from the proof above that, for each $K>0$,
\bestar
S_{nK} &:=&\frac 1{\sqrt n}\sum_{j=1}^n H_K(x_{nj}, u_{nj})\, \ep_{nj}\, \to_D \,\int_0^1 \widetilde H_K(X_t)dB_{t}+\si_K(X)\, {\cal N} ,
\eestar
where $ \si_K^2(X)=\int_0^1 \big[\widehat H_K^2(X_t)- \widetilde H_K^2(X_t)\big] dt  $  and ${\cal N}$ is a standard normal variate independent of $B=\{B_t\}_{t\ge 0}$ and  $X=\{X_t\}_{t\ge 0}$. This yields (\ref {m1}) since
$
\mathbb{P}(S_n\not =S_{nK}) \le  \mathbb{P}\Big(\max\limits_{1\le k\le n} |x_{nk}|> K\Big) \to 0
$
and
\bestar
&&\mathbb{P}\Big(\int_0^1 \widetilde H(X_t)dB_{t}+\si(X)\, {\cal N}\not=\int_0^1 \widetilde H_K(X_t)dB_{t}+\si_K(X)\, {\cal N} \Big)\no\\
&\le & \mathbb{P}\Big(\sup_{0\le t\le 1} |X_t|> K\Big) \to 0,
\eestar
 as $n\to\infty$ first and then $K\to\infty$. The proof of Theorem \ref {main} is complete.
\hfill $\Box$

\subsection {Proof of Theorem \ref {th8} }\la {sec4.2}

We start with two lemmas  on the asymptotics of general nonlinear LS estimators. Consider a nonlinear regression model having the form:
\be  \la{model3}
y_k=g_k(\theta) +u_k, \quad k=1,2, ..., n,
\ee
where  $\theta=(\theta_1, ...,\theta_q)$ is a vector of  unknown parameters, $\{u_k, {\mathcal F}_k\}_{k\ge 1}$ is a martingale difference with $E(u_k^{2}\mid {\cal F}_{k-1})<\infty, a.s.$  and  $g_k(\theta)$ for each $\theta\in \Theta$, a compact of ${\mathbb R}^q$,  is adapted to $ {\cal F}_{k-1}$.
Let  $\widehat \theta_n$  be the LSE of  $\theta$ and $\theta_0$ be   an interior of $\Theta$. The following consistency and asymptotic normality results come from   \cite{JW1}. See also \cite{W3}.

\begin{lem} \la {cort}  Suppose that  there exist a sequence of   constants $0<k_n\to\infty$ and a sequence of random variables  $T_k$  (adapted to $ {\cal F}_{k-1}$)   such that
\begin{enumerate}
\item[(a)] $|g_k(\theta_1)-g_k(\theta_2)|\le ||\theta_1-\theta_2||^\al\, T_k$   for all $\theta_1, \theta_2\in \Theta$ and for some $\al>0$;
\item[(b)]  $\frac 1{k_n}\sum_{k=1}^n  T_k^2\,  =O_P(1)$ and
$\sum_{k=1}^n\, T_k^2 \big[1+E(u_k^2\mid {\cal F}_{k-1})\big]  =O_P(k_n^2/\log^2 n)$;
\item[(c)] the finite dimensional distributions of $\frac 1{k_n}\sum_{k=1}^n \,\big[g_k(\theta)-g_k(\theta_0)\big]^2$ converge to those of  $ G(\theta)$, where $G(\theta), \theta\in \Theta$ is a stochastic process  of $\theta$ satisfying
 $P(\inf_{||\theta-\theta_0||\ge \de\atop \theta\in \Theta}G(\theta)>0)=1$   for each $\de>0$.
\end{enumerate}
Then
$
||\widehat \theta_n-\theta_0||=o_P(1),
$
i.e.,  $\hat\theta_n$ is a consistent estimator of $\theta_0$.
\end{lem}

 Let $\dot{g}_k(\theta)= \big(\frac {\partial g_k(\theta)}{\partial \theta_1}, ...,\frac {\partial g_k(\theta)}{\partial \theta_q}\big)^{\prime}$  be the first  derivative of $g_k(\theta )$, $Z_{n}(\theta)=(D_{n}^{-1})^{\prime }\,\sum_{k=1}^{n}\, \dot{g}_k(\theta)\,u_k$ and
  \bestar
  Y_{n}=(D_{n}^{-1})^{\prime }\,\sum_{k=1}^{n}\, \dot{g}_k(\theta_0)
\dot{g}_k(\theta_0)^{\prime }\,D_{n}^{-1},
  \eestar
  where $D_n=\mbox{diag}(d_{1n},...,d_{qn})$ is a sequence of diagonal matrices satisfying  $ n^{-\delta}\,\min_{1\le j\le q} d_{jn}\to \infty$ for some $\delta>0$.

\begin{lem} \la {cor1?} Suppose that
\begin{itemize}
\item[(i)] $||D_n^{-1} \big[\dot{g}_k(\theta_1)-\dot{g}_k(\theta_2)\big]|| \le  ||\theta_1-\theta_2||^{\al}\, T_{nk}$  for some $0<\al\le 1$ and  for any $\theta_1, \theta_2\in \Theta$,  where $T_{nk}$ is   adapted to  ${\mathcal F}_{k-1}$ for each $n\ge 1$, satisfying
\be
\sum_{k=1}^n\, T_{nk}^2\big[1 +E(u_k^2 | {\cal F}_{k-1})\big]+\sum_{k=1}^n \,T_{nk}^2=O_P(1); \la {A3.21}
\ee
\item[(ii)] $\widehat \theta_n\to_P \theta_0$,  $Z_{n}(\theta_0)=O_P(1)$ and $Y_{n}\rightarrow _{D}M$, where $
M>0$, a.s., i.e., the smallest eigenvalue of $M$ is almost surely positive.
\end{itemize}
We have
\be\la{thm3.1.1}
D_{n}(\hat{\theta}_{n}-\theta _{0})=Y_{n}^{-1}\,Z_{n}(\theta_0)+o_{P}(1).  \label{thm3.11}
\ee
If in addition $(Y_{n}, Z_{n})\rightarrow _{D} (M, Z)$
where $M>0$, a.s., then   $D_{n}(
\hat{\theta}_{n}-\theta _{0})\rightarrow _{D}M^{-1}\,Z$.

\end{lem}

We turn back to the proof of Theorem \ref {th8}  by checking the conditions of   Lemma  \ref {cor1?} with
 $$ g_k(\theta) := f(x_k, \Delta x_{k}, ..., \Delta x_{k-d}, \  \theta)=f(x_k,  \xi_{k}, ..., \xi _{k-d}, \  \theta). $$

We first show that $\widehat \theta_n\to_P \theta_0$ by using Lemma \ref {cort}  with $k_n=nv_f^2(d_n)$. In this case,
it only needs to show  the following results:  for any $\theta_j\in \Theta$ and $\al_j\in R$, $j=1,2,..., q$,
\be
\sum_{j=1}^q\al_j\, \frac 1{nv_f^2(d_n)}\sum_{k=1}^n \big[g_k(\theta_j)-g_k(\theta_0)\big]^2& \to_D& \sum_{j=1}^q\al_j\,G(\theta_j),  \la {e11}
\ee
where $G(\theta):=
\int_0^1 \big[ \widetilde h_f(X_t, \theta)-\widetilde h_f(X_t, \theta_0)\big] dt$ and
\be
&&\frac 1{n}\sum_{k=1}^n \big[1+(|x_k|/d_n)^{2\ga}\big] (1+||\Gamma_k||^{2\beta})(1+\si_k^2) =O_P(1), \la {e21}
\ee
where $\Gamma_k=(\xi_k,...,\xi_{k-d})$.
Indeed, due to the continuity of $X_t$, it follows from  (\ref {tc1}) that\\  $ P\big[\inf_{||\theta-\theta_0||\ge \de\atop \theta\in \Theta}G(\theta)>0\big]=1$.
 This, together with (\ref {e11}), implies the condition (c) of Lemma \ref {cort}.
On the other hand,     the conditions (a) and (b) of Lemma \ref {cort} follow easily from Assumption \ref {as3.3}[(i) and (iii)] with $T_k=T_f(x_k)(1+||\Gamma_k||^\beta)$  and (\ref {e21}).  We therefore have $\widehat \theta_n\to_P \theta_0$.

We only prove   (\ref {e11}), where the idea is similar to that of \cite{W3}. Result (\ref {e21}) can be obtained by (\ref {ad35}) with $H(x, y)=(1+|x|^{2\ga})(1+||y||^{2\beta})$ and  hence the details are omitted.
Write $x_t^*=x_tI(|x_t|/d_n\le A)$,
\bestar
g_k^*(\theta) &=& f(x_k^*, \xi_k,...,\xi_{k-d}, \  \theta), \no\\
R_n(\theta) &=& \sum_{k=1}^n \big[g_k( \theta)-v_f(d_n)h_f(x_k/d_n, \xi_k,...,\xi_{k-d}, \theta)\big]^2 ,\no\\
R_n^*(\theta) &=& \sum_{t=1}^n \big[g_k^*(\theta)-v_f(d_n)h_f(x_t^*/d_n, \xi_k,...,\xi_{k-d},\theta)\big]^2.
\eestar
 For any fixed $A>0$, it follows from  Assumption \ref{as3.3}(ii) and (iii) that, as $n\to\infty$,
\bestar
\sup_{\theta\in \Theta}|R_n^*(\theta)|&\le&
o(1)\, v_f^2(d_n) \sum_{t=1}^nT_{f}^2(x_t^*/d_n)(1+||\Gamma_k||^{2\beta})\no\\
 &=&o(1) \sum_{t=1}^n(1+|x_t^*/d_n|^{2\ga})(1+||\Gamma_k||^{2\beta}) =o_P(1)\, n v_f^2(d_n).
\eestar
 This implies that, for any $\ep>0$,
\bestar
&& P\big(\frac 1{n v_f^2(d_n)}\sup_{\theta\in \Theta}|R_n(\theta)| \ge \ep \big) \no\\
&\le &
 P \big(x_k\not=x_k^*, \mbox{for some $$k=1,...,n} \big)+ P
\big(\frac 1{n v_f^2(d_n)}\sup_{\theta\in \Theta}|R_n^*(\theta)|\ge \ep\big)  \no\\
&\le & P\big(\max_{1\le k\le n}|x_k|/ d_n\ge A\big) +  P
\big(\frac 1{n v_f^2(d_n)}\sup_{\theta\in \Theta}|R_n^*(\theta)|\ge \ep\big)  \no\\
&\to& 0,
\eestar
as $n\to\infty$ first and then $A\to \infty$, namely, we have
\be
\sup_{\theta\in \Theta}|R_n(\theta)| &=&
o_P\big[n v_f^2(d_n)\big]. \la {e32}
\ee
Now, by letting
\bestar
\Delta_n &=&\sup_{\theta\in \Theta}|R_n(\theta)| +\sup_{\theta\in \Theta}|R_n(\theta)|^{1/2}\,\no\\
&&\qquad
\Big(\frac 1{n}\sum_{k=1}^n \big[h_f(x_k/d_n, \xi_k,...,\xi_{k-d}, \theta)-h_f(x_k/d_n, \xi_k,...,\xi_{k-d}, \theta_0)\big]^2\Big)^{1/2},
\eestar
and noting that, for any $\theta\in \Theta$,
\bestar
&& \frac 1{nv_f^2(d_n)}\sum_{k=1}^n \big[g_k (\theta)-g_k(\theta_0)\big]^2\no\\
&= &
\frac 1{n}\sum_{k=1}^n \big[h_f(x_k/d_n, \xi_k,...,\xi_{k-d}, \theta)-h_f(x_k/d_n,\xi_k,...,\xi_{k-d}, \theta_0)\big]^2 +4\,\widetilde\Delta_n,
\eestar
where $|\widetilde\Delta_n|\le \Delta_n$,
result (\ref {e11}) follows from (\ref {e32}) and (\ref {ad35}). 

We now turn back to the proof of Theorem \ref {th8} by  verifying other conditions of Lemma  \ref {cor1?}. First of all, by letting
\bestar
T_{nk}=\frac 1{\sqrt n}\sum_{j=1}^qv_{\dot {f}_j}^{-1}(d_n) T_{\dot {f}_j} (x_k) (1+||u_k||^{\beta}),
\eestar
it follows easily from Assumption 3.3 (iii) and (\ref {e21})  that
\bestar
&& \sum_{k=1}^n\, T_{nk}^2\big[1 +E(u_k^2 | {\cal F}_{k-1})\big]+\sum_{k=1}^n \,T_{nk}^2 \no\\
&\le& \frac 1{n}\sum_{t=1}^n \big[1+(|x_t|/d_n)^{2\ga}\big] (1+||u_t||^{2\beta})(1+\si_t^2) =O_P(1).
\eestar
This yields (\ref {A3.21}).  So, together with  $\widehat \theta_n\to_P \theta_0$ that is proved above,    it remains to show that, for any $\al_i'=(\al_{i1},...,\al_{iq})\in \mathbb{R}, i=1,2,3, $
\be
&& (\alpha _{1}^{\prime }\,Y_{n}\alpha _{2},\ \ \alpha _{3}^{\prime }Z_{n})\no\\
&&\rightarrow _{D}\ \Big( \alpha _{1}^{\prime }\int_0^1\,\Psi(t)\Psi(t)^{\prime}\,dt\,\alpha _{2},\ \alpha _{3}^{\prime }\,Z\,\Big), \la {e31}
\ee
where  $Z:=\int_0^1 \widetilde G_t\,dB_t+ {\Omega}^{1/2}\, {\cal N}_q$ and, with $D_n=\sqrt n\, \mbox{diag}\, \big(v_{\dot {f}_1}(d_n), ...,v_{\dot {f}_q}(d_n)\big)$,
\bestar
   Y_{n}=(D_{n}^{-1})^{\prime }\,\sum_{k=1}^{n}\, \dot{g}_k(\theta_0)
\dot{g}_k(\theta_0)^{\prime }\,D_{n}^{-1}, \quad Z_{n}=(D_{n}^{-1})^{\prime }\,\sum_{k=1}^{n}\, \dot{g}_k(\theta_0)\,\si_k\ep_k, .
  \eestar

Write $ h_{\dot f}(x, y, \theta)=\big( h_{\dot{f}_1}(x,y;\theta),\,..., \,h_{\dot{f}_q}(x,y;\theta)\big)$ and
\bestar
 \widetilde  Y_{n} &=& \frac 1n \,\sum_{k=1}^{n}\,\alpha _{1}^{\prime }\,h_ {\dot f}(x_k/d_n, \xi_k,...,\xi_{k-d}, \theta_0)h_{\dot f}(x_k/d_n, \xi_k,...,\xi_{k-d}, \theta_0)^{\prime }\,\alpha _{2}\,, \\
 \widetilde Z_{n} &=& \frac 1{\sqrt n}\,\sum_{k=1}^{n}\, \alpha _{3}^{\prime }\,h_{\dot f}(x_k/d_n, \xi_k,...,\xi_{k-d}, \theta_0)\,\si_k\ep_k .
  \eestar
Note that, by (\ref {ad35}) in Remark \ref {r2.2},
\bestar
 ( \widetilde  Y_{n},\ \ \widetilde Z_{n} )
&\rightarrow _{D}& \Big( \alpha _{1}^{\prime }\int_0^1\,\Psi(t)\Psi(t)^{\prime}\,dt\,\alpha _{2}, \ \alpha _{3}^{\prime }\,Z\,\Big).
\eestar
and
\bestar
\widetilde  Y_{n} &=& \frac 1n \,\sum_{k=1}^{n}\,\alpha _{1}^{\prime }\,h_ {\dot f}(x_k/d_n, \xi_k,...,\xi_{k-d}, \theta_0)h_{\dot f}(x_k/d_n, \xi_k,...,\xi_{k-d}, \theta_0)^{\prime }\,\alpha _{2}\,, \\
 \widetilde Z_{n} &=& \frac 1{\sqrt n}\,\sum_{k=1}^{n}\, \alpha _{3}^{\prime }\,h_{\dot f}(x_k/d_n, \xi_k,...,\xi_{k-d}, \theta_0)\,\si_k\ep_k .
\eestar
Result (\ref {e31}) will follow if we prove: for all $i, j=1,..., q,$
\be
R_{1n}(i, j)& :=& \frac 1{ n v_{\dot f_i}(d_n)v_{\dot f_j}(d_n)}  \sum_{k=1}^n \big[\dot{f}_i(...)\dot{f}_j(...)-v_{\dot f_i}(d_n)v_{\dot f_j}(d_n)h_{\dot{f}_i}(...)h_{\dot{f}_j}(...)\big] \no\\
&=& o_P(1) , \la{jk9} \\
R_{2n}(j) &:=& \frac 1{\sqrt n v_{\dot f_j}(d_n)}  \sum_{k=1}^n \big[\dot{f}_j(...)-v_{\dot f_j}(d_n)h_{\dot{f}_j}(...)\big]\si_k \, \ep_k\ =\ o_P(1). \la {jk10}
\ee

We only prove (\ref {jk10}). The proof  of (\ref {jk9}) is similar to that of  (\ref {e32}) and hence the details are omitted. To this end, for any $A>0$,  let
\bestar
R_{2n}^*(j) &=& \frac 1{\sqrt n v_{\dot f_j}(d_n)}  \sum_{k=1}^n \big[\dot{f}_j(...)-v_{\dot f_j}(d_n)h_{\dot{f}_j}(...)\big]I(|x_k|/d_n\le A)\si_k \, \ep_k.
\eestar
Note that, as $n\to\infty$ first and then $A\to \infty$,
\bestar
P(R_{2n}(j) \not = R_{2n}^*(j) ) &\le &P \big(\max_{1\le k\le n}|x_k|/d_n\ge A \big) \to 0.
\eestar
It suffices to show that, for each $A>0$, $R_{2n}^*(j) =o_P(1).$ This follows easily from  the following fact: in term of   Assumption \ref {as3.3} (ii) and (iii),
\bestar
\E \big[R_{2n}^*(j)\big]^2 &=&  \frac 1{ n v_{\dot f_j}^2(d_n)}  \sum_{k=1}^n \E \big\{\big[\dot{f}_j(...)-v_{\dot f_j}(d_n)h_{\dot{f}_j}(...)\big]^2I(|x_k|/d_n\le A)\si_k^2\big\} \no\\
&\le &\sup_{\theta\in \Theta, y\in R}\frac {|R_{\dot f_j}(d_n\, x, y,  \theta)|}{T_{\dot f_j}(d_n x)(1+||y||^\beta)}\,
\frac {1+A}n \sum_{k=1}^n\E\big[ (1+||\Gamma_k||^{2\beta})\si_k^2\big] \to 0,
\eestar
as $n\to \infty$, where $\Gamma_k=(\xi_k,...,\xi_{k-d})$. The proof of Theorem \ref {th8} is complete.  \hfill $\Box$

\subsection{Proof of Lemma \ref {lem31}} \la {subs1}
 Let $l_{nj}=g_N^2(u_{nj})$. Recall that $H(x, y)=0$ for all  $|x|\ge N$, implying  $\widetilde H(x)$ is bounded  (by $C_0$, say). It follows from Assumption \ref {as2.2}(b) that, for any $K>0$,
\be
 \max_{1\le k\le n} |y_{nj}|
&\le & \frac {1}{\sqrt n}\max_{1\le k\le n}\big[\,g_N(u_{nj})\,|\ep_{nj}| I_{(|l_{nj}|\ge K)}\big]+  \frac {C_0+K}{\sqrt n} \max_{1\le j\le n} |\ep_{nj}| \no\\
 &:=&  R_{1n} +(C_0+K)\, R_{2n}. \la {4.6ap}
\ee
Using the uniformity of $l_{nj}$ and $E(\ep_{nj}^2|{\cal F}_{n, j-1})=1$, we have
\bestar
\E  R_{1n} &\le& \(\frac 1n \sum_{j=1}^n \E \big[g_N^2(u_{nj})\, \ep_{nj}^2 I_{|l_{nj}|\ge K}\big] \)^{1/2} \no\\
&= &\( \frac 1 n \sum_{j=1}^n \E  l_{nj}^2I_{|l_{nj}|\ge K} \)^{1/2}\to 0,
\eestar
as $n\to\infty$ first and then $K\to \infty$, i.e, $ R_{1n}=o_P(1)$. On the other hand,  it follows from $ \frac 1{\sqrt n }\sum_{j=1}^{[nt]} \ep_{nj} \Rightarrow B_t$ on $D_R[0, \infty)$ that
\bestar
R_{2n} &=&\frac {1}{\sqrt n} \max_{1\le j\le n} |\ep_{nj}| =o_P(1).
\eestar
See, Proposition VI 3.26 of Jacod and Shiryaev (2001), for instance.  Taking these estimates into (\ref {4.6ap}), we prove $\max_{1\le j\le n} |y_{nj}| =o_P(1)$, as required in  (\ref {4.6a}).
The proof of Lemma \ref  {lem31} is now complete.  \hfill $\Box$

\subsection{Proof of Lemma \ref {lem3}} \la {subs1a}
 We start with some preliminaries. Since  the limit process $X_t$ is path continuous,  the convergence $x_{n,[nt]}\Rightarrow X_t$ on $D_{\mathbb{R}}[0,1]$ can be understood in the sense of the
uniform topology.  It follows from this fact   that
\begin{eqnarray}
\limsup_{N\rightarrow\infty}\limsup_{n\rightarrow\infty}\Big[\,\mathbb{P}\,\big(%
\max\limits_{1\le j\le n} |x_{nj}|\ge N\big)+\mathbb{P}\,\big(%
\max\limits_{1\le j\le n} |x_{nj}|\le 1/N\big)\Big]=0.  \label{leK}
\end{eqnarray}
Furthermore, by the tightness of $\{x_{n,[nt]}\}_{ 0\le t\le 1}$, for any $%
\varepsilon>0$ and $\delta>0$, there is some $\tilde{\delta}=\tilde{\delta}%
(\varepsilon,\delta)>0$ such that
\begin{eqnarray}
\mathbb{P}\,\Big(\,\sup_{|s-t|\le \tilde{\delta}} |x_{n,[nt]}-x_{n,[ns]}| \ge \delta\,\Big)\le
\varepsilon  \label{leK1}
\end{eqnarray}
 for all sufficiently large $n$. In terms of (\ref{leK1}), for any $%
\delta>0$, we have
\begin{eqnarray}
\lim_{n\to\infty }\mathbb{P}\,\Big(\,\max_{0\le k\le n/m} \max_{km\le j\le (k+1)m}
|x_{nj}-x_{n,km}|\ge \delta\,\Big)= 0,  \label{atightness2}
\end{eqnarray}
for any $m:=m_n \to \infty$ satisfying $n/m\to\infty$. On the other hand, it follows from (\ref {c1.2}) and the continuities of $H(x, y)$ and $\widetilde H(x)$ (with respect to $x$) that, as $m\to \infty$,
\be
\max_{m\le k\le n-m} \E \sup_{|x| \le N}\, \Big|\frac 1m \sum_{j=km+1}^{(k+1)m}\, H(x, u_{nj})-\widetilde H(x)\Big|\to 0, \la {c2.20}
\ee
for any fixed $N>0$.

We are now ready to prove (\ref {4.19}).
Let $H_1(x, y)=H(x, y)-\widetilde H (x)$. Note that $H(x, y)$ is uniformly  continuous  on  $|x|\le N$  and $H_1(x, y)=0$ for all $|x|\ge N$ and $y\in \mathbb{R}$. It follows   that, for each $\ep>0$,  there exists an $\de_{\ep}$ depending only on $N, \ep$ only such that
\be
|H_{1}(x_1, y)-H_{1}(x_2, y)|&\le& |H(x_1, y)-H(x_2, y)| \,
\le \ep g_N(y), \la {5.201}
\ee
whenever $|x_1-x_2|\le \de_\ep$, where $g_N(y)=\sup_{|x|\le N} H(x, y)$.
 Write $$\Omega_{\delta_\ep}=\Big\{\omega: \max\limits_{0\le k\le n/m} \max\limits_{km\le j\le (k+1)m}  |x_{nj}-x_{n,km}|\le \delta_{\ep}\Big\}$$ and $T_n=[n/m]-1$, where $m$ is chosen so that $m\to\infty$ and $n/m\to\infty$.
In terms of  (\ref {5.201}), on $\Omega_{\delta_\ep}$,  we have
\bestar
\widehat R_{n, 1} &:=&\frac{1}{n}\, \max\limits_{0\le k\le n/m}\,  \Big|\sum_{j=1}^{(k+1)m}\, H_{1}(x_{nj}, u_{nj})\, \Big|
\ \le \  \frac{1}{n} \sum_{k=0}^{T_n}
\Big| \sum_{j=km+1}^{(k+1)m}\, H_{1}(x_{nj}, u_{nj})\,\Big| 
\no\\
&\le &  \frac{1}{n} \sum_{k=0}^{T_n}
\Big| \sum_{j=km+1}^{(k+1)m}\, H_{1}(x_{n, km}, u_{nj})\,\Big| 
\no\\
&&\qquad +\,
\frac{1}{n}\,\sum_{k=0}^{T_n}
\, \sum_{j=km+1}^{(k+1)m}\,\big|\, H_{1}(x_{nj}, u_{nj})-H_{1}(x_{n, km}, u_{nj})\,\big| \, \no\\
&\le&  \frac{1}{n}\,\sum_{k=0}^{T_n}  \sup_{|x|\le N}
\Big|\frac 1 m \sum_{j=km+1}^{(k+1)m}\, H_{1}(x, u_{nj})\,\Big|   + \frac {\ep}{n}\sum_{j=1}^{n} g_N (u_{nj}).
\eestar
Therefore,  for any $\eta, \eta_1>0$ and  $\ep=\eta\, \eta_1$, it follows from  Assumption \ref {as2.2}(b), (\ref {atightness2}) and  (\ref {c2.20}) that
\bestar
\mathbb{P}(|\widehat R_{n, 1}|\ge \eta) &\le&\mathbb{P}\big (\bar  \Omega_{\delta_\ep}\big) +\mathbb{P}(|\widehat R_{n, 1}|\ge \eta,  \Omega_{\delta_\ep}) \no\\
&\le& \mathbb{P}\big (\bar  \Omega_{\delta_\ep}\big) + \eta^{-1}\,\max_{m\le k\le n-m} \E  \sup_{|x|\le 2N}
\Big|\frac 1 m \sum_{j=km+1}^{(k+1)m}\, H_{1N}(x, u_{nj})\,\Big| \no\\
&& +
\frac {\ep\, \eta^{-1}}{n}\sum_{j=1}^{n} \mathbb{E} g_N (u_{nj}) \ \le\ C\,\eta_1 +o(1),
\eestar
 where $\bar  \Omega_{\delta_\ep}$ denotes the complementary set of $  \Omega_{\delta_\ep}$.
 This proves  $\widehat R_{1n, 1}=o_P(1)$.

Similarly, by recalling $|H_{1}(x, y)|\le  g_N (u_{nj})$, it follows from Assumption \ref {as2.2}(b) that
\bestar
\widehat R_{n, 2} &:=&\frac{1}{n}\, \max\limits_{0\le k\le n/m}\, \sum_{j=km+1}^{(k+1)m}\,\big| H_{1}(x_{nj}, u_{nj})\, \big|\no\\
 &\le&\, \frac {1}n \max_{0\le k\le n/m} \, \sum_{j=km+1}^{(k+1)m} g_N (u_{nj}) \no\\
 & \le& \frac {m}n \max_{1\le j\le n } g_N (u_{nj})=o_P(1).
\eestar
In terms of  these estimates, we have
\bestar
&& \frac 1n \max_{1\le k\le n} \big| \sum_{j=1}^{k}\big[H(x_{nj}, u_{nj}) -\widetilde H(x_{nj})\big]\big| \no\\
&\le& \frac{1}{n}\, \max\limits_{0\le k\le n/m}\, \max\limits_{km\le j\le (k+1)m} \Big|\sum_{k=1}^{j}\, H_{1}(x_{nj}, u_{nj})\, \Big|\no\\
&\le& \frac{1}{n}\, \max\limits_{0\le k\le n/m}\,  \Big|\sum_{j=1}^{(k+1)m}\, H_{1}(x_{nj}, u_{nj})\, \Big|+\frac{1}{n}\, \max\limits_{0\le k\le n/m}\, \sum_{j=km+1}^{(k+1)m}\,\big| H_{1}(x_{nj}, u_{nj})\, \big|
\no\\
&=& \widehat R_{n, 1} +\widehat R_{n, 2}=o_P(1),
\eestar
establishing  (\ref {4.19}).

We next prove (\ref {3.4}).
Since $K_j^2(x, y)=H^2(x, y)-2H(x, y) \widetilde H(x)+ \big[\widetilde H(x)\big]^2$, we may write
\bestar
V_{nk}^2- \widetilde V_{nk}^2 &=&\frac 1n\sum_{j=1}^k\big[  H^2(x_{nj}, u_{nj})^2-\widehat H^2(x_{nj})\big] +\frac 2n\sum_{j=1}^k\widetilde H (x_{nj})\big[\widetilde H (x_{nj})-H(x_{nj}, u_{nj}) \big]\no\\
&:=& A_{nk}+B_{nk}.
\eestar
The similar argument as in the proof of (\ref {4.19}), but $H(x, y)$ is replaced by $H^2(x, y)$ and $\widetilde H(x) \,H(x, y)$ respectively, yields that
\bestar
\max_{1\le k\le n}\big|V_{nk}^2-\widetilde V_{nk}^2\big|&\le&\max_{1\le k\le n}\big| A_{nk}\big|+\max_{1\le k\le n}\big| B_{nk}\big| =
o_P(1)
\eestar
as required.
\hfill $\Box$

\subsection {Proof of Lemma \ref {lem4}} \la {subs3}
 We first prove (\ref {1.12y}). First note that, for each $s\ge 1$,
 \bestar
 \frac 1n \sum_{j=n+1}^{[ns]}l_j^2 \to_{a.s.} s-1.
 \eestar
 This implies that $\tau_n(t)\le n$ if $0<t\le V_{nn}^2$ and, when $t> V_{nn}^2$,
 \bestar
 n\le \tau_n(t) \le n(1+2t), \quad a.s.
 \eestar
 In term of these facts, for each $t>0$, we have
  \be
 \max_{1\le j\le \tau_n(t)} |\widehat y_{nj}| &\le & \max_{1\le j\le n} | y_{nj}| +\max_{n\le j\le n(1+2t)} | l_j|/\sqrt n=o_P(1), \la {fin1}
 \ee
 wehre we have used (\ref {4.6a} ) and the fact that, due to $\E|l_j|^3\le C<\infty$,
 \bestar
 \max_{n\le j\le n(1+2t)} | l_j|/\sqrt n \le \Big(\frac 1{n^{3/2} }\sum_{j=n}^{n(1+2t)}|l_j|^3\big)^{1/3}=o_P(1).
 \eestar
 Similarly, for each $t>0$, it follow from Assumption \ref {as2.2}(b)  that
\be
\E \big|\sum_{j=1}^{\tau_n(t)}\E\big(\widehat y_{nj}^2|{\cal F}_{n, j-1}\big)-t\big| &\le& \frac 1n \E\max_{1\le j\le n}
|V_{nj}^2-V_{n, j-1}^2|+ \frac 1n \E \max_{n\le j\le n(1+2t)}l_j^2 \no\\
&\le &\frac Cn \E \max_{1\le j\le n}  g_N^2 (u_{nj})\,   +  \frac 1n \E \max_{n\le j\le n(1+2t)}l_j^2\no\\
&\to &  0. \la {fin2}
\ee
As a consequence,   we have  $\sum_{j=1}^{\tau_n(t)} \widehat y_{nj}\Rightarrow B'_t $ on $D_{R}[0, \infty)$ by the classical martingale limit theorem.

We finally prove (\ref {1.12a}).
 It is readily seen that  $\{Z_{n}(t)\}_{t\ge 0}$ is tight. Therefore, for each $\{n'\}\subseteq \{n\}$, there exists a subsequence    $\{n''\}\subseteq \{n'\}$ such that, on $D_{\mathbb{R}^4}[0, \infty)$,
\bestar
Z_{n''}(t) &\Rightarrow& Z_t:= \big( B_t,\  B_{1t},\ X_t, B_t'\big).
\eestar
As a consequence, (\ref {1.12a}) will follow if we prove that $B'$ is independent of $X$  and $F$.  Note that $\big(\,\ep_{nj}, \,\eta_{nj}, \, \widehat y_{nj}, \, {\cal F}_{nj}\big)_{n\ge 1, j\ge 1}$ forms a  martingale array at the same probability space. The independence between $B'$ and ${\mathrm F}$ follows the standard martingale limit theorem and the fact: by using  (\ref {4.19}) and definition of $\widehat y_{nj}$,
\bestar
 && \frac 1 {\sqrt n} \max_{1\le k\le \tau_n(t)} \big|\big| \sum_{j=1}^{k} \E \big(\xi_{nj} \widehat y_{nj} |{\cal F}_{n, j-1}\big)\big|\big|\no\\
 &=&\frac 1 {\sqrt n} \max_{1\le k\le n} \big|\big| \sum_{j=1}^{k} \E \big(\xi_{nj}  y_{nj} |{\cal F}_{n, j-1}\big)\big|\big| \to 0,
 \eestar
 where $\xi_{nj}=\big(\,\ep_{nj}, \,\eta_{nj}\big)$. The fact that $B'$ is independent of $X$ is obvious since $X=\{X_t\}_{t\ge 0}$ is   adapted to the filtration ${\mathrm F}=\{\si (B_s, B_{1s}, s\le t)\}_{t\ge 0}$. The proof of Lemma \ref {lem4} is complete.
\hfill $\Box$

\bibliographystyle{ecta}
\bibliography{REF-martingale-2026}

\end{document}